\documentclass[12pt]{article}
\usepackage{amsmath}
\usepackage{amssymb}
\usepackage[english,francais]{babel}

\newtheorem{theorem}{Theorem}[subsection]
\newtheorem{lemma}[theorem]{Lemma}
\newtheorem{e-proposition}[theorem]{Proposition}

\newtheorem{e-definition}[theorem]{Definition\rm}

\newtheorem{propo}{Proposition}[subsubsection]

\title{Error structures and parameter estimation}
\author{Nicolas Bouleau, Christophe Chorro}
\date{Nov 2004}
\begin{document}

\maketitle

\noindent{\bf Abstract}
This article proposes and studies a link between statistics and the theory of Dirichlet forms used to compute errors. The error calculus based on  Dirichlet forms is an extension of  classical Gauss' approach to error propagation.  The aim of this paper is to derive error structures from measurements. The links with Fisher's information lay the foundations of a strong connection with experiment. Here we show that this connection behaves well towards   changes of variables and is related to  the theory of asymptotic statistics. Finally the study  of products  permits to lay the premise of an infinite dimensional empirical error calculus.

\selectlanguage{english}

\textbf{Mathematical subject classification (2000):} 31C25, 47B25, 49Q12, 62F99, 62B10, 65G99.\\

\textbf{Keywords:} Error, sensitivity, Dirichlet forms, squared field operator, Cramer-Rao inequality, Fisher information.
\section{\textbf{Introduction}}
\subsection{\textbf{Intuitive notion of error  structures}}
Let us consider a random quantity  $C$ (for example the concentration of some pollutant in a river) that can be measured by an experimental device which result exhibits an error denoted by $\bigtriangleup C$. These quantities may be represented as random variables generally correlated (for higher pollution levels, the device becomes fuzzier). In this classical probabilistic approach we have to know the law of the pair $(C,\bigtriangleup C)$ or equivalently the law of  $C$ and the conditional law of  $\bigtriangleup C$  given $C$. Thus, the study of error transmission is associated to the calculus of images of probability measures. Unfortunately, the knowledge of the law of  $\bigtriangleup C$ given $C$  by means of experiment is practically impossible. Now, let us look at the propagation of errors when the errors are small. For the sake of simplicity we adopt temporarily the following assumptions:\\
$\bullet $ Only  the conditional variance $var[\bigtriangleup C\mid C]$ is known.\\
$\bullet$ The errors are small enough to allow the simplification usually performed by physicists: $\bigtriangleup C=\varepsilon Y$ where  $Y$ is a bounded random variable and  $\varepsilon$  a size parameter.\\

 If  $f$ is   $C^3(\mathbb{R},\mathbb{R})$  with bounded derivatives, supposing at first that the error is conditionally centered ${E}[\bigtriangleup C\mid C]=0$, Taylor's formula gives $$\bigtriangleup (f(C))=f'(C)\bigtriangleup C + \frac{1}{2} f^{''}(C) (\bigtriangleup C)^2 +\varepsilon^3~ 0(1)  $$ hence \begin{center}$var[\bigtriangleup f(C)\mid C ]=f'^2(C) var[\bigtriangleup C\mid C] + \varepsilon^3~ 0(1) $\end{center}\vspace{0.2cm}\begin{center}$~\mathbb{E}[\bigtriangleup f(C)\mid C]=\frac{1}{2} f''(C)var[\bigtriangleup f(C)\mid C] + \varepsilon^3~ 0(1).$\end{center}~

 In the same way, for another regular function $ h$  we have :

  \begin{equation}var[\bigtriangleup (h\circ f(C))\mid C ]=h'^2(f(C)) var[\bigtriangleup f(C)\mid C] + \varepsilon^3~ 0(1)  \end{equation}\begin{equation}\mathbb{E}[\bigtriangleup (h\circ f(C))\mid C]= h'(f(C)) \mathbb{E}[ \bigtriangleup f(C)\mid C] + \frac{1}{2} h''(f(C))var[\bigtriangleup f(C)\mid C] + \varepsilon^3~ 0(1). \end{equation}

These formulae of the  propagation of variances and biases show that once a nonlinear function has been applied, the error is no longer centered and the bias has the same order of magnitude as the variance. Through other applications this phenomenon persists. Moreover we can see that the calculus on the variances is a first order calculus and does not involve the biases whereas the calculus on the biases is of  second order and involves the variances. This remark is fundamental:    the error calculus on variances is necessarily the first step of an analysis of errors based on differential methods. It will be the main focus of our study.\\

 On the probability space  associated to the observation of $C$, $(\mathbb{R},Bor(\mathbb{R}), law~of~C)$ ( where $Bor(\mathbb{R})$ is the borelian $\sigma$-field of $\mathbb{R}$),  we introduce the  operator  $\Gamma^C$ called the quadratic error operator which provides for each function  $f$ the asymptotical conditional variance of the error on  $f(C)$: $$\Gamma^C[f](x)=\underset{{\tiny \varepsilon\rightarrow 0}}{\lim}\frac{var[\bigtriangleup f(C)\mid C=x]}{\varepsilon^2}.$$ As the covariance operator in probability theory, $\Gamma^C$ polarizes into a bilinear operator: $$\Gamma^C[f,g](x)= \underset{{\tiny \varepsilon\rightarrow 0}}{\lim}\frac{covar[\bigtriangleup f(C),\bigtriangleup g(C)\mid C=x]}{\varepsilon^2}.$$ Moreover if F is in  $C^2(\mathbb{R}^2,\mathbb{R})$ with bounded derivatives,  we obtain a transport formula known as the Gauss' law of errors propagation (\cite{bib:5}, Chap.1, Appendix): \begin{equation}\Gamma^C[F(f,g)]=F_1'^2(f,g) \Gamma^C[f]+F_2'^2(f,g) \Gamma^C[g]+2F'_1(f,g)F'_2(f,g) \Gamma^C[f,g].\end{equation}

Now we can adopt an intuitive definition of an error structure:

An error structure is a probability space  $(W,\mathcal{W},m)$ equipped with a positive, symmetric, bilinear  operator  $\Gamma$ acting on random variables and  fulfilling a first order functional calculus on regular functions: $$\Gamma[F(f_1,\ldots,f_n)]=\sum\limits_{i,j} F_i'(f_1,\ldots,f_n)F_j'(f_1,\ldots,f_n)\Gamma[f_i,f_j].$$

If  $\phi:\mathbb{R}\rightarrow \mathbb{R}$ is a regular mapping, this definition is preserved by image: we can equip the image space $(\mathbb{R},Bor(\mathbb{R})$, law of $\phi(C))$ with the quadratic error operator $\Gamma^{\phi(C)}$ associated to the observation of $\phi(C)$. We have the following fundamental relation \begin{equation}\Gamma^{\phi(C)}[f](x)=\mathbb{E}[\Gamma^C[f(\phi)](C)\mid \phi(C)=x].\end{equation}

When we observe a two-dimensional quantity $C=(C_1,C_2)$ with   erroneous components  modelled with two error structures   $(\mathbb{R},Bor(\mathbb{R})$, law of $ C_1,\Gamma^{C_1})$ and  $(\mathbb{R},Bor(\mathbb{R})$, law of $C_2$,$\Gamma^{C_2})$,  if  $(C_1,\bigtriangleup C_1)$ is independent of $(C_2,\bigtriangleup C_2)$ we need to define an error structure $(\mathbb{R}^2,\mathbb{B}(\mathbb{R}^2)$, law of $C_1$ $\otimes $ law of $C_2$, $\Gamma^{{\tiny C_1\otimes C_2}})$ such that   $\Gamma^{{\tiny C_1\otimes C_2}}$ expresses a summation of errors component per component. Indeed, if   $F:\mathbb{R}^2\rightarrow \mathbb{R}$ is regular, from the independence hypothesis it follows    $$var[\bigtriangleup (F(C_1,C_2))\mid (C_1,C_2)]=F_1'^2(C_1,C_2)var[\bigtriangleup C_1\mid C_1] + F_2'^2(C_1,C_2)var[\bigtriangleup C_2\mid C_2] + \varepsilon^3 0(1)$$ thus \begin{equation}\Gamma^{C_1\otimes C_2}[F](x,y)=\Gamma^{C_1}[F(.,y)]+\Gamma^{C_2}[F(x,.)].\end{equation}

The preceding intuitive considerations lead to the following rigorous mathematical framework.

\subsection{\textbf{An extension tool}}
Now we present an  axiomatic extension of the preceding notion of error structures using the language of Dirichlet forms. It gives a powerful  tool easy to handle in error calculations and sensitivity analysis. As noticed above, we limit ourselves to a first order calculus which is already significant in most of applications. We refer to  \cite{bib:5} for a calculus on biases involving the infinitesimal generator associated to the underlying Dirichlet form. This error calculus based on Dirichlet forms lies between the probabilistic approach (errors are supposed to be random variables) and the deterministic one (dealing with infinitely small deterministic errors to use differential calculus).\\

From now on, an error structure is a term $(W,\mathcal{W},m,\mathbb{D},\Gamma)$ where $(W,\mathcal{W},m)$ is a probability space, $\mathbb{D}$ is a dense vector subspace of  $L^2(m)$ and  $\Gamma$ is a positive symmetric bilinear map from     $\mathbb{D}\times\mathbb{D}$ into $L^1(m)$ fulfilling:
\begin{itemize}
\item[1)] the functional calculus of class  $C^1 \cap Lip$ i.e. if  $U=(U_1,\ldots,U_n)\in \mathbb{D}^ n $, $V=(V_1,\ldots,V_p)\in \mathbb{D}^p $,  $F\in C^1(\mathbb{R}^n,\mathbb{R})\cap Lip=\{C^1$ and Lipschitz$\}$ and $G\in C^1(\mathbb{R}^p,\mathbb{R})\cap Lip$ then, 
$$(F(U_1,\ldots,U_n),G(V_1,\ldots,V_p))\in \mathbb{D}^2$$ and
$$\Gamma[F(U_1,\ldots,U_n),G(V_1,\ldots,V_p)]=\sum\limits_{i,j} F_i'(U) G_j'(V)\Gamma[U_i,V_j],$$\item[2)] $1\in\mathbb{D}$ (this implies $\Gamma[1]=0$), \item[3)] the bilinear form     $\mathcal{E}[F,G]=\frac{1}{2}\int \Gamma[F,G] dm $ defined on $\mathbb{D}\times \mathbb{D}$ is closed i.e. $\mathbb{D}$  is complete under the norm of the graph $$\parallel. \parallel_{\mathcal{E}}=\\(\parallel.\parallel^2_{L^2(m)}+\mathcal{E}[.])^{\frac{1}{2}}.$$ We always write $\Gamma[F]$ for $\Gamma[F,F]$ and $\mathcal{E}[F]$ for $\mathcal{E}[F,F]$.
\end{itemize}~\\

This  notion  is derived from the theory of Dirichlet forms (\cite{bib:2} Ch.1,\cite{bib:8},\cite{bib:**}). It is a natural extension of the classical Gauss approach (\cite{bib:3}) and  it seems to be a good way to study the propagation of errors and the sensitivity to changes of parameters in physical and financial models (\cite{bib:3},\cite{bib:4},\cite{bib:5}).

The condition 1) is similar to the Gauss' law of  small errors propagation  (3). For $U=(U_1,\ldots,U_n)\in \mathbb{D}^ n$, the intuitive meaning of the matrix   $\underset{=}{\Gamma} [U]=[\Gamma[U_i,U_j]]_{1\leq i,j\leq n}$ is the variance-covariance of  the error on  $U$ (\cite{bib:5} Ch.1). Implicitly, we still suppose that the error is infinitely small although it is not mentioned  in the notation. It is as if we had an infinitely small unit to measure errors that was fixed in the whole problem. Then, the   hypothesis   $3)$ is added to the heuristic definition and can be seen as a coherence principle. In fact, if the random variables  $(X_n)_{n\in\mathbb{N}}$ and $X$ are in  $\mathbb{D}$, if  $X_n\rightarrow X$ in $L^2(m)$ and  ($X_n$, error on $X_n$) converges in a suitable sense, it converges necessarily  to the pair ($X$, error on $X$).

From the hypotheses mentioned above,  $\mathcal{E}$ is a local Dirichlet form and $\Gamma$  its associated squared field operator. The domain  $\mathbb{D}$ is preserved by Lipschitz functions: if  $F:\mathbb{R}^n\rightarrow \mathbb{R}$ is a contraction in the following sense $$|F(x)-F(y)|\leq\sum\limits_{i=1}^{n}|x_i-y_i|$$ then for  $U=(U_1,\ldots,U_n)\in \mathbb{D}^n$ one has $F(U)\in \mathbb{D}$ and $$\Gamma[F(U_1,\ldots,U_n)]^{\frac{1}{2}}\leq \sum\limits_{i=1}^{n}\Gamma[U_i]^{\frac{1}{2}}.$$
We would like to emphasize that the closedness property is the key stone of our approach. It plays the same role as the $\sigma$-additivity in probability theory and permits to compute the errors on  functions known as limits of simpler objects.

The operations of taking  images by mapping (definition 3.1.2)     and making  countable products (definition 5.0.8)  naturally provide error structures on spaces of stochastic processes (\cite{bib:2} Ch.2,\cite{bib:4},\cite{bib:5} Ch.6).

Since a probability space $(W,\mathcal{W},m)$ can be known thanks to statistical experiments, we raise the problem of the empirical identification of an error structure. In the same way as the $\sigma$-additivity of $m$ on $\mathcal{W}$ could not result from experiments but is a fundamental mathematical hypothesis, our error structure will have to verify the  closedness property 3) (This cannot be deduced from observation). Thus let $\theta$ be a parameter taking its values in an open set $\Theta\subset\mathbb{R}^d$. It is frequently useful to treat $\theta$ as the realization of a random variable $V:(\Omega,\mathcal{A},\mathbb{P})\rightarrow \Theta$ with a known distribution $\rho$ chosen by combining experience with convenience (\cite{bib:12} p.225). Let $X$ be a random variable defined on the probability  space $(\Omega,\mathcal{A},\mathbb{P})$ with values in a measurable space $(E,\mathcal{F})$. Let us denote by    $P_\theta$ the conditional law of $X$ given $V=\theta$. Classically,  to estimate $\theta$ we may use the statistical model $(P_\theta)_{\theta\in\Theta}$ generated by the observations of   $X$.   Here we want to equip $\Theta$ with an error structure  \begin{equation}S^{V}=(
\Theta, \mathcal{B}(\Theta), \rho   ,
\mathbb{D}^{V},    \Gamma^{V})\end{equation} where $\Gamma^{V}$ will express  the precision of our knowledge on $\theta$. Our approach is to consider $\Gamma^{V}$ as the inverse of the Fisher matrix  which is an accuracy measure  for  regular statistical models (see \cite{bib:7}). We will study the behavior  of this identification through changes of variables and products to show its remarkable stability.

\section{The Cramer-Rao Inequality (C.R.I.) and the Fundamental Identification (F.I.)}

\subsection{\textbf{Regular models.}}

From now on  $(.,.)$ will denote the usual scalar product on $\mathbb{R}^d$ and $\parallel.\parallel$ its associated norm. We suppose that  $(P_{\theta})_{\theta\in \Theta}$  satisfies the  conditions of  regular models (\cite{bib:10} p.65):\\ \begin{itemize}
\item[(a)] The measures $P_{\theta}$  are absolutely continuous with respect to a $\sigma$-finite measure $\mu$ and $\frac{dP_\theta}{d\mu}=f(.,\theta)>0$.
\item[(b)]  $\theta \rightarrow {f(x,\theta)}$ is continuous for   $\mu$-almost all $x$.
\item[(c)] We set  $g(x,\theta)=\sqrt{f(x,\theta)}$.  There exists $\phi:E\times \Theta\rightarrow \mathbb{R}^d$ such that$~\forall~\theta\in\Theta$,
$$\int \parallel \phi(x,\theta)\parallel ^2 d\mu(x)<\infty$$ and
$$\int |g(x,\theta+h)-g(x,\theta)-(\phi(x,\theta),h)|^2 d\mu(x) =o(\parallel h \parallel^2).$$ thus  the positive semi-definite  matrix $J(\theta)= 4\int \phi(x,\theta)\phi(x,\theta)^t d\mu(x)$ is defined as  the Fisher information matrix of our model.
\item[(d)] $\theta\rightarrow \phi(.,\theta)$ is continuous in  $L^2(\mu)$.
\item[(e)] The model is  identifiable:  $\theta\rightarrow P_{\theta}$ is injective.
\end{itemize}~\\

\textbf{Remarks A:} \textbf{i)} There exists several definitions of regular models. Here we use a notion taken from  \cite{bib:10}  where the  conditions are quite general. These hypotheses are made to allow a differentiation under integrals which is needed for the proof of the Cramer-Rao inequality. We can found in  \cite{bib:1} another definition using the classical differential calculus and supposing that $J$ is continuous when it is a simple consequence of d).

\textbf{ii)} The assumption  c) is a condition of differentiability in quadratic mean in  $L^2(\mu)$.  Moreover, if we assume that  $\theta\rightarrow f(.,\theta)$ is differentiable in the classical sense then $\phi(x,\theta)=\frac{\nabla
f(x,\theta)}{2\sqrt{f(x,\theta)}}$ and we obtain the following   expression of the so-called Fisher information matrix
$$J(\theta)=\left [\int \frac{\frac{\partial f(x,\theta)}{\partial \theta_i} \frac{\partial f(x,\theta)}{\partial \theta_j}}{f(x,\theta)} d\mu(x)\right ]_{0\leq i,j\leq d}.$$ To establish the differentiability in quadratic mean, one often proceeds by showing classical differentiability and equi-integrability  (see  \cite{bib:6},\cite{bib:13}).

 \textbf{iii)}  Identifiability   is a purely statistical hypothesis. Intuitively, it means  that the model  can distinguish two different values of the parameter  $\theta'\not=\theta''$ if and only if  $P_{\theta'}\not=P_{\theta''}$. In this case,  if independent experiments are available,  we have an infinite  family of independent variables with the same law  $P_\theta$ denoted by $Z^{\theta}=(X^{\theta}_i)_{i\in\mathbb{N}}$ and for $\theta'\not=\theta''$, the laws of the processes $Z^{\theta'}$ and $Z^{\theta''}$ are mutually singular. Thus,  $\theta'$ and  $\theta''$ are perfectly identified thanks to experiment.\hfill$\Box$ \\

\subsection{\textbf{Cramer-Rao Inequality}}

\begin{theorem}$ \mbox{\rm{(\cite{bib:10} p.73)}} $  Let $\psi:\mathbb{R}^d\rightarrow \mathbb{R}^m$ be differentiable and  $(P_{\theta})_{\theta\in\Theta}$ be a regular model with $\forall\theta\in\Theta$ $det (J(\theta))\not=0$. If  $T(X)$ is an unbiased estimator of $\psi(\theta)$ such that $\mathbb{E}[T(X)^2\mid V=\theta]$ is locally bounded in $\theta$   then
$$\mathbb{E} [(T(X)-\psi(\theta))(T(X)-\psi(\theta))^t\mid V=\theta]\geq \psi'(\theta)J^{-1}(\theta)\psi'(\theta)^t.$$
where $\geq$  is the order relation between symmetric matrices defined by the cone of positive symmetric ones.
\end{theorem}~\\

\textbf{Remark B:} An estimator  $T(X)$ fulfilling the hypotheses of the preceding theorem is said to be a regular unbiased estimator of  $\psi(\theta)$.\hfill$\Box$\\

Now, up to  the end, we suppose  that the Fisher information  matrix is regular. Thus, the Cramer-Rao inequality gives a bound of estimation for the quadratic risk. Let us have a look on the error structure $(6)$ we want to determine. If the components of identity are in $\mathbb{D}^{V}$, according to the functional calculus, we have  for $F\in {Lip}^{1}(\Theta)=\{ F \in C^1(\Theta,\mathbb{R})$ and Lipschitz$\}$,   $$ \Gamma^V[F]=(\nabla F)^t \underset{=}{\Gamma}^{V}[Id] (\nabla F)$$ where the matrix  $ \underset{=}{\Gamma}^{V}[Id] (\theta)$ represents the error of estimation on  $V$ given $V=\theta$. Since $ \underset{=}{\Gamma}^{V}$ takes its significance from a calculus on variances, the Cramer-Rao inequality leads us to state the  fundamental identification
 \begin{equation*}
\underset{=}{\Gamma}^{V}[Id] =J^{-1}.~~~~\quad\bf{(F.I)}
\end{equation*}
 As well as the statistical identification of a probability space presupposes the  $\sigma$-additivity of the measure, we want to determine an error structure deriving from experiment in which  $\mathcal{E}^{V}$ is a closed form. According to the fundamental identification we make the following assumption:\\

\textbf{Hypothesis} (\textbf{E}): From now on, we suppose  the existence of a dense vector subspace of $L^2(\rho)$ denoted by $\mathbb{D}^V$ and the existence of an operator $\Gamma^V$ fulfilling conditions $1)$, $2)$ and $3)$  such that  $Lip^1(\Theta) \subset\mathbb{D}^{V} $ and,    for all $ F$ in $Lip^1(\Theta)$,  $\Gamma^{V}[F]=F'J^{-1}(F')^t$. Moreover, as    $\mathbb{D}^{V} $  may not be uniquely defined, we take it minimal for inclusion, which implies the density of  ${Lip^1(\Theta)}$ in $ \mathbb{D}^{V}$  for the norm $\parallel.\parallel_{\mathcal{E}^{V}}$.\\

This hypothesis dictates conditions on $\rho$ and $J^{-1}$ which are often  fulfilled as seen in the following proposition  (see also \cite{bib:8}):\\

\begin{e-proposition}
a)  Let  $\Theta$  be a bounded open set of $\mathbb{R}^d$ of the form $\Theta=\prod\limits_{i=1}^d ]\theta_0^i,\theta_1^i[$  where  the $\theta_j^i$ are real numbers such that $\theta_1^i\geq\theta_0^i$. We shall assume that  $\rho$ is a probability measure which is absolutely continuous with respect to the  Lebesgue measure on  $\Theta$ with a positive  density $q$ in $Lip^1(\Theta)$. Suppose that the model $(P_\theta)_{\theta\in \Theta}$ can be extended to a regular model on an open set  $\Theta'$ such that $\overline{\Theta}\subset\Theta'$. Then,  hypothesis \textbf{E} is fulfilled.\\

b) When $\Theta=\mathbb{R}$, if we assume that   $\int_\Theta \theta^2 d\rho(\theta)<\infty$ and that $\frac{1}{J}$ belongs to  $L^1(\rho)$, the hypothesis \textbf{E} is equivalent to the conditions of Hamza theorem (\cite{bib:8} p.105).

\end{e-proposition}~\\

\textbf{Proof:} a) Let $(F_n)_{\mathbb{N}}$ be a sequence in  $Lip^1(\Theta)$ such that  $F_n\rightarrow 0$ in $L^2(\rho)$ and  $\Gamma^V[F_n-F_m]\rightarrow 0 $ in $L^1(\rho)$ where  $\Gamma^V:  Lip^1(\Theta)\rightarrow L^1(\rho) $ is well-defined by  $\Gamma^V[F]=F'J^{-1}(F')^t$. If we show that $\Gamma^V[F_n]\rightarrow 0$ in  $L^1(\rho)$, the conclusion follows according to  \cite{bib:8} p.4.

One  defines the mapping $\Phi$: $$\Phi:\left( \begin{tabular}{l}
$ \overline{\Theta} \times \mathcal{S}_d \longmapsto \mathbb{R}^*_{+} $\\ $(\theta,\xi) \longrightarrow \sum\limits_{i,j=1}^{d}a_{i,j}(\theta)\xi_i\xi_j $
\end{tabular}\right)  $$ where  $\mathcal{S}_d$ is the unit sphere of $\mathbb{R}^d$ and where the coefficients of $J^{-1}$ are denoted by $a_{i,j}$ .

The function $\Phi$ is continuous on a compact set, thus there exists $\delta,\delta'>0$ such that $\delta\leq\Phi\leq\delta'$. It implies
$$\delta  |\nabla F_n -\nabla F_m|^2     \leq\Gamma^V[F_n-F_m]\leq \delta' |\nabla F_n -\nabla F_m|^2 .~~\qquad(*)$$
Hence,  $\nabla F_n$ is a cauchy sequence in $L^2(\rho;\mathbb{R}^d)$ and there is a function  $G=(G_1,\ldots,G_d)$  in $L^2(\rho;\mathbb{R}^d)$ satisfying for all $i\in\{1,\ldots,n\}$ $\partial_i F_n\rightarrow G_i$ in $L^2(\rho)$.\\
Let $ \phi$ be a function in $C^{\infty}_K(\Theta)=\{ F \in C^\infty(\Theta,\mathbb{R})~\mbox{\rm{with compact support}}\}$. One notices that $\phi,q,F_n$ are Lipschitz and can be extended to $\overline{\Theta}$. Thus,  by integration by parts formula we obtain  $$-\int_\Theta F_n (\partial_i\phi) q d\theta=\int_\Theta (\partial_iF_n) \phi q d\theta+\int_\Theta F_n \phi (\partial_i q) d\theta $$

 and by passing to the limit,  it follows that $\forall \phi\in C^{\infty}_K(\Theta),$  $$\int_\Theta G_i \phi q d\theta=0$$ thus $G_i=0$.  We can conclude using the inequality $ (*).$   \\

b) Hamza theorem gives necessary and sufficient conditions for the existence of an error structure $S=(\mathbb{R}, \mathcal{B}(\mathbb{R}), \rho,
\mathbb{D},\Gamma)$  such that  $C^{\infty}_K(\mathbb{R})\subset\mathbb{D}$ and $\Gamma[F]=\frac{F'^2}{J}$ on $C^{\infty}_K(\mathbb{R})$.

Let  $(F_n)_{n\in\mathbb{N}}$ be a sequence in  $C^{\infty}_K(\mathbb{R})$ with the same Lipschitz constant $1$ such that $F_n\rightarrow Id$ everywhere with $\forall n $ $|F_n|\leq|Id|$ and  $F_n'\rightarrow 1$ everywhere. Using the dominated convergence theorem and  the closedness of $\Gamma$  we obtain that $Id\in\mathbb{D}$, hence $Lip^1(\mathbb{R})\subset\mathbb{D}$ and $\Gamma[F]=\frac{F'^2}{J}$ for  $F\in Lip^1(\mathbb{R})$. The result follows naturally.  \hfill$\Box$\\

\textbf{Remarks C:} \textbf{i)}  The statistical situation with a constant information matrix is often encountered in  classical parametric models (see \cite{bib:12}):  Location family, Normal models with fixed coefficient of variation, Logistic model, Scale parameter. In this case the condition of extension of the model $(P_\theta)_{\theta\in \Theta}$ can be removed in a).

 \textbf{ii)}  The operator $\Gamma^{V}$ is bilinear. It is possible to introduce a new operator, the gradient, denoted by $\nabla^{V}$, which can be seen as a signed and linear version of the standard deviation of the error  and satisfies $\forall F\in\mathbb{D}^{V}$ $$\Gamma^{V}[F]=(\nabla^{V}[F],\nabla^{V}[F]).$$

Since the error structure $S^{V}$ is defined on a finite dimensional space it is easy to construct  $\nabla^{V}$ putting $$\nabla^{V}:\left( \begin{tabular}{l}
$\mathbb{D}^{V} \longmapsto L^2(\rho;\mathbb{R}^d)$ \\ $ F \longrightarrow R (F')^t$
\end{tabular}\right)  $$where $R$ is the square root of  $J^{-1}$. The gradient fulfills the classical differentiation  chain rule.

\textbf{iii)} We can notice that the fundamental identification gives, without other hypotheses, a second order calculus with variances and biases as mentioned in the introduction. In fact, we can associate to the  Dirichlet form $\mathcal{E}^{V}$ a unique self adjoint operator $A^{V}$ (see  \cite{bib:2},\cite{bib:8}), called the infinitesimal generator. It has a domain  $D(A^{V})$ included in $\mathbb{D}^{V}$ and  it takes its values in  $L^1(\rho)$.  Moreover we have $$A^{V}[F(U)]=F'(U)A^{V}[U]+\frac{1}{2}F''(U)\Gamma^{V}[U]$$ when $U\in D(A^V)$, $\Gamma^{V}[U]\in L^2(\rho)$ and $F:\Theta\rightarrow \mathbb{R}$  is a function of class  $C^2$ with bounded derivatives. Thus, the preceding formula expresses the propagation of the conditional expectation of the error  in the same way as (2).\hfill$\Box$\\

Now, we want  to test the robustness of  the fundamental identification by comparing its  properties with the well-known behavior of  the Fisher information in the classical framework of parametric estimation.

\section{\textbf{Change of variables: the injective case.} }

We are going to show the stability of the fundamental identification for regular changes of variables.
\subsection{\textbf{The regular injective case.}}
\begin{e-definition}We suppose that  $\psi:\Theta\rightarrow \mathbb{R}^d$ is injective of class $ \mathcal{C}^{1}\cap Lip$. This change of variables is said to be  regular if  $det (\psi'(x))\not=0$    for all $x$.\\
\end{e-definition}
From the local inversion theorem, it follows  that $\psi$  is a  $\mathcal{C}^1$-diffeomorphism  on its image and  $\psi(\Theta)$ is an open set of $\mathbb{R}^d$.\\
Now, we want to equip $\psi(\Theta)$ with an error structure that expresses the intrinsic accuracy of our knowledge on $\psi(\theta)$. There are two natural ways to proceed.

\subsubsection{\textbf{From the estimation point of view. }}~\\
In the injective case, the change of variables is just a reparameterisation of the model. To estimate  $\psi(\theta)$ we use the model $(P_{\psi^{-1}(a)},~a\in \psi(\Theta)
)$. Since $dP_{\psi^{-1}(a)}(x)=f(x,\psi^{-1}(a)) d\mu(x)$, we can see easily that this model is regular. Let us have a look on  the error structure we obtain using the fundamental identification. The operator $ \Gamma^{\psi(V)}$  is defined on $
Lip^{1}(\psi(\Theta))$ by
$$ \Gamma^{\psi(V)} [{F}](a)=(\nabla_a{F})^t (J^{\psi(V)}(a))^{-1} (\nabla_a {F})~~~~\forall a \in \psi(\Theta) ~~ $$
 where $J^{\psi(V)}$ is the Fisher information matrix of the regular model $(P_{\psi^{-1}(a),~a\in \psi(\Theta) })$. Moreover, as  $\forall a \in \psi(\Theta)$, $$J^{\psi(V)}(a)=
 [\psi'(\psi^{-1}(a))^{-1} ]^t  ~ [J(\psi^{-1}(a) )]~[\psi'(\psi^{-1}(a))^{-1} ]$$   one has for
 $F \in {Lip}^{1}(\psi(\Theta))$
$$ \Gamma^{\psi(V)} [{F}](a)=(\nabla_a{F})^t  [\psi'(\psi^{-1}(a)) ]  ~ [J(\psi^{-1}(a) )]^{-1}~[\psi'(\psi^{-1}(a)) ]^t (\nabla_a {F}).$$ Using that  $\psi$ is injective of class  $C^1\cap Lip$, from hypothesis \textbf{E} it follows that the form $\mathcal{E}^{\psi(V)}$ defined  on $Lip^1(\psi(\Theta))$ by $$\mathcal{E}^{\psi(V)}[F]=\frac{1}{2} \int \Gamma^{\psi(V)}[F] d\psi_*\rho$$ is closable and we denote by  $\mathbb{D}^{\psi(V)}$ the domain of its smallest closed extension. Thus, the error structure associated to the fundamental identification for the estimation of  $\psi(\theta)$ is $$S^{\psi(V)}=(\psi(\Theta),
\mathcal{B}(\psi(\Theta)), \psi_*\rho,
\mathbb{D}^{\psi(V)},
 \Gamma^{\psi(V)}) .$$

\textbf{Remark D:}  When $d=1$ one obtains
$$J^{\psi(V)}(\psi)=  \frac{J}{\psi'^2}.$$

Hence, if  $\psi$ is flat enough in  $\theta$,  $\psi(\theta)$ can be estimated more accurately than $\theta$. This property is intuitively coherent since a value $\theta'$ at a given small distance from $\theta$ will lead to a smaller deviation of $\psi(\theta')$ from $\psi(\theta)$ the smaller the value of $|\psi'(\theta)|$ is.\hfill$\Box$
 \subsubsection{\textbf{From the error calculus point of view}}~\\
Among the advantages of the error calculus based on Dirichlet forms, let us emphasize here its practical  flexibility. It is easy to define both the product of error structures and the image of an error structure by a mapping. The following definition is the rigorous formulation of the intuitive expression  (4) which corresponded to a change of observation in our preliminary   study of error calculus.\\

\begin{e-definition}

Let $S=(W, \mathcal {W}, m,
{\mathbb{D}},\Gamma)$ be an error structure and    $Y:W\rightarrow
\mathbb{R}^d$ $\in$ $\mathbb{D}^d$ such that $Y(\mathcal{W})$ is an open set of $\mathbb{R}^d$. Let us define $ \widetilde{\mathbb{D}_{Y}}= \lbrace f\in {L}^{2}( Y_*m)\mid f(Y)\in  \mathbb{D} \rbrace$ and for  $f\in\widetilde{\mathbb{D}_{Y}}$,
$\widetilde{\Gamma_{Y}}[f](x)=\mathbb{E}_m[\Gamma[f(Y)]\mid  Y=x]$.

 If we denote by $\mathbb{D}_Y$  the closure of $Lip^1(Y(\mathcal{W}))$ in $(\widetilde{\mathbb{D}_Y},\parallel.\parallel_{\widetilde{\mathcal{E}_{Y}}}  )$ and by   $\Gamma_Y$ the restriction of $\widetilde{\Gamma_{Y}}$ to ${\mathbb{D}_Y}$ then $${\psi_*S}=(Y(\mathcal{W}), \mathcal{B}(Y(\mathcal{W})), Y_*m,
{\mathbb{D}_{Y}}, {\Gamma_{Y}})$$ is an error structure called the image structure of $S$ by $Y$.
\end{e-definition}~

Let us study the image of $S^{V}$ by $\psi$ which is another natural way to endow $\psi(\Theta)$ with an error structure. For  $F\in {Lip}^{1}(\psi(\Theta))$ one has, $\forall a\in Im(\psi)$,
$$\mathbb{E}_{\tiny{\rho}}[\Gamma^{V}[F(\psi)]\mid ~\psi
=a]=\mathbb{E}_{\tiny{\rho}}[   \nabla({F}(\psi))^t ~
J^{-1}~ \nabla({F}(\psi))\mid ~\psi=a]~~~~$$ and

$${\Gamma^{\psi(V)}}[F]= \Gamma_{\psi}^{V}[F ]~~~~\psi_*\rho~~a.e.$$

Thus, ${\Gamma^{\psi(V)}}$ and
${\Gamma_{\psi}^{V}}$ are equal on  $ Lip^{1}(\psi(\Theta))$.

Using the density of $Lip^1(\psi(\Theta))$ in $(
\mathbb{D}^{V}_{\psi},\parallel.\parallel_{\mathcal{E}^{V}_{\psi}})$, we have the following expected property:\\
\begin{e-proposition}~The fundamental identification is preserved by the transformation $\psi$. In other terms:
$$\psi_*S^{V}=S^{\psi(V)}.$$ \end{e-proposition}

\textbf{Remark E:} Suppose we are studying the sensitivity of a physical or financial model depending on the parameter $\theta$ to small random perturbation by using an error structure on $\Theta$ and the functional calculus for $\Gamma$ to compute the propagation of errors on the outputs of the model. If the error structure is obtained by the Fisher information matrix of a statistical model as above, the preceding invariance result means that the accuracy on $\theta$ has a physical significance, independently of mathematical repameterization.\hfill$\Box$

\subsection{\textbf{The non-regular injective case}}

After the regular case studied in the preceding section, let us  see what happens at a point $\theta$ such that $\psi'(\theta)$ is singular. First, we supposes that $d=1.$

Let $a_{0}$ be  equal to $\psi(\theta_{0})$ with $\theta_0 \in \Theta$
and $\psi'(\theta_0)=0$. We can see easily that the model   $(P_{\psi^{-1}(a)},~a\in \psi(\Theta) )$ possesses an irregularity at $a_0$. Intuitively, as far as estimation is concerned, this situation is not harmful because it induces a good approximation of $a_0$ (see Remark D). If we put  $J^{\psi(V)}(a_0)= +\infty$ it follows
$$\Gamma^{\psi(V)}(Id)(a_0)=\frac{1}{J^{\psi(V)}(a_0)}=0=\Gamma_{V}^{\psi}(Id)(a_0).$$

In the general  case, since $J(\theta_0)$ is supposed to be definite positive, we can reduce simultaneously $\psi'(\theta_0)$ and $J(\theta_0)$ and work component per component.
If $\psi'(\theta_0)$ is singular, there exists eigenvectors for to the eigenvalue  $0$ which correspond to directions of infinite information for $J^{\psi(V)}(a_0)$. The other  eigendirections are dealt as in the regular case.

 We can see that the fundamental identification is  still stable  in this case.\\

\textbf{Remarks F:} \textbf{i)} The concept of infinite information  appears in asymptotic statistics where it expresses a faster convergence of the maximum likelihood estimator toward the parameter.

\textbf{ii)} We have seen that, for injective changes of variables, the error structure obtained estimating directly
 $\psi(\theta)$ coincides  with the image by  $\psi $ of the structure associated to the estimation of $\theta$. This phenomenon can be viewed as a sufficiency principle (well known for the Fisher information \cite{bib:10}, p.70) because when $\psi$ is injective, $P_{\theta}$ depends on  $\theta$ only through  $\psi$.

\textbf{iii)} The proposition $3.1.3$ is based on the simple relation between  $J^{\psi(V)}$ and  $J$. This property of the Fisher information is not fulfilled  for other types of information bound. For example, for the bounds of  Bhattacharaya type (see \cite{bib:10}) (which involve higher derivatives and are more precise) it is impossible to obtain such a coherence property. The hypotheses of regular models are the good level of axiomatization for our study.\hfill$\Box$

\section{\textbf{The non injective case.}}

We are now in a special situation: we have put in correspondence an error structure and a parametric model thanks to the Fisher information. But on one side (error structures) non-injective changes of variables are allowed (def 3.1.2) and on the other side (statistical models) they meet difficulties. We derive benefit of this remark to propose a new framework  for the estimation of a parameter in this case which is directly linked with the notion of error structure.

Here we suppose that  $\psi$ is a   function in $Lip^1(\Theta)$  not necessarily injective but such that $\psi(\Theta)$ is an open set of $\mathbb{R}^d$ (in order to apply definition $3.1.2$ with $Y=\psi$ ).

To estimate  $\psi(\theta)$ the reparameterisation introduced in the previous section is meaningless. To avoid this problem, we give a new protocol.
\subsection{\textbf{Estimation protocol of  $\psi(\theta)$ when  $\psi$ is not injective.}}

 To estimate  $\theta$ we use a regular model  $(P_{\theta})_{\theta\in \Theta}$ such that hypothesis \textbf{E} is fulfilled. In this section the random variables  $X_1\ldots X_n$ defined on $(\Omega,\mathcal{A},\mathbb{P})$ with values in $(E,\mathcal{F})$ will be, given  $V=\theta$, a n-sample of $P_\theta$.

 To estimate  $\psi(\theta)$, it is natural to  use the model   $(Q_a)_{a\in
\psi(\Theta)}$ generated by the observation of $X_1\ldots X_n$ given  $\psi(V)=a$. From the definition of the conditional expectation it follows that
$$dQ_a(x)=\mathbb{E}_{\rho}[f(x,.)\mid \psi=a] d\mu (x).$$In particular, we need a global knowledge of  $f(x,.)$   to perform  $Q_a$.~\\

\textbf{Remarks G:}   \textbf{i)} When $\psi$ is injective, the preceding protocol coincides with the reparameterisation. In this case a pointwise knowledge is sufficient.

\textbf{ii)} The case of  non-injective  changes of parameters is often tackled  in the literature on the following restrictive form. When  d=1, we consider a point  $\theta_0$ such that  $\psi'(\theta_0)\not=0$.
Since  $\psi$ is in $C^1(\Theta,\mathbb{R})$, according to the local inversion theorem there exists  $\theta^{min}_0$, $\theta^{max}_0$ such that
$\psi: \small{]\theta^{min}_0, \theta^{max}_0[ \rightarrow
\psi( ]\theta^{min}_0, \theta^{max}_0[) }$ \normalsize is a
$C^1$-diffeomorphism with an inverse denoted by $i_{\theta_0}$ (which depends on  $\theta_0$ contrary to the injective case).

If we suppose that  previous observation leads us to believe that $\theta$ is in
$]\theta^{min}_0, \theta^{max}_0[$,  locally we are going back to the case processed in section $3$. Thus, we set  $\forall a \in
\psi( ]\theta^{min}_0, \theta^{max}_0[)$,  $$J^{\psi(V)}_{\theta_0}(a)=
\frac{J(i_{\theta_0}(a))}{\psi'^2(i_{\theta_0}(a))}    .$$

The quantity $J^{\psi(V)}_{\theta_0}(a)$ is called the local Fisher information because it takes into account only one antecedent of  $a$.

When we do not have any \textit{a priori} information on $\theta$, one has to use a concept which expresses the entire behavior of  $\psi$.\hfill$\Box$\\     %The operator $\Gamma^{V}_{\psi}$ The On doit donc en quelque sorte effectuer un %passage du local au global qui fera intervenir de manière naturelle l'opérateur $\Gamma^{V}_{\psi}$.             % qui ne  fera pas %intervenir l'Information de Fisher ( qui revet un caractère intrinséquement local) mais  le calcul %d'erreur.

  Since  $\psi$ is non-injective, the model $(Q_a)_{a\in\psi{(\Theta)}}$ may present irregularities. Thus, the Fisher information may be undefined. Moreover, even if it exists, the information matrix is not easy to perform. So we are going to show the relevance of error calculus in this case, showing that the operator $\Gamma^{V}_\psi$ is a substitute of the inverse of the Fisher information in the sense that it gives a simple bound of estimation and is linked to asymptotic statistics.

\subsection{ \textbf{$\Gamma_{\psi}^{V}$ as an estimation bound.}}

To simplify, let us suppose that $d=1$.

Using a regular parametric model to estimate  $\theta$, we have seen that for a regular unbiased estimator  $T(X)$ of  $\psi(\theta)$, the Cramer-Rao inequality
\begin{equation} \mathbb{E}[ (T(X)-\psi(\theta))^2\mid V=\theta]\geq \frac{\psi'^2(\theta)}{J(\theta)}\end{equation}  gave a bound of the quadratic risk and lead to interpret  $J$ as the information  on $\theta$ contained in  observation  $X$.
 In the same way, when the estimators are built with  the independent observations  $(X_1,\ldots X_n)$, it is easy to see that the additivity property  of the Fisher information matrix ensures that \begin{equation*}\mathbb{E}[
(T(X_1,\ldots,X_n)-\psi(\theta))^2\mid V=\theta]\geq
\frac{\psi'^2(\theta)}{n J(\theta)}\end{equation*} if  $\mathbb{E}[T(X_1,\ldots X_n)|V=\theta]=\psi(\theta)$.

Thus, conditioning  (7) with respect to  $\psi$ one has
$$\mathbb{E}[ (T(X)-a)^2\mid \psi(V)=a]\geq \mathbb{E}_{\rho}[\frac{\psi'^2}{J}\mid \psi=a]=\Gamma_{\psi}^{V}[Id](a)$$and  ${\Gamma^{V}_{\psi}[Id]}$ appears as a natural bound of the problem. Similarly, one obtains$$\mathbb{E}[
(T(X_1,\ldots,X_n)-a)^2\mid \psi(V)=a]\geq
\frac{\Gamma^{V}_{\psi}[Id](a)}{n} .$$

\textbf{Remark H:}  $\frac{1}{\Gamma^{V}_{\psi}}$  can be seen as an  additive information  when independent observations are combined.\hfill$\Box$\\

\subsection{\textbf{Links with asymptotic statistics.}}
For the sake of simplicity about the question of existence and  unicity of the  maximum likelihood estimator we suppose that for the model $(P_{\theta})_{\theta\in\Theta}$, for all $n\in\mathbb{N}$, for all  $(x_1,\ldots,x_n)\in E^n$, the equation $$ \sum\limits_{i=1}^{n} \frac{\partial}{\partial \theta} logf(x_i,\theta)=0$$
 has a unique solution denoted by  $\widehat{\theta_n}(x_1,\ldots,x_n)$ which is a maximum for the function  $\theta\rightarrow \prod\limits_{i=1}^{n} f(x_i,\theta)$. In this section we assume that  $\Theta$ is a convex  bounded subset of $\mathbb{R}$ (this could easily be extended to any finite dimension).

  In order to show that $\Gamma_{\psi}^V$ is the key stone of some asymptotic results, one requires preliminary knowledge concerning the convergence of  the sequence of estimators $(\widehat{\theta_n}(X_1,\ldots,X_n))_{n\in\mathbb{N}}$.
\subsubsection{\textbf{Convergence of the maximum likelihood estimator.}}~

We essentially refer the reader to \cite{bib:9}, \cite{bib:10} for the proof of the results exposed here
and for complementary details.

The asymptotic techniques used in this section can be easily  extended to a more general framework than the case of experiments based on the observation of n-samples (especially for the applications to stochastic processes). These techniques are not based on the historical approach using Taylor's formula  any more (see for  example  \cite{bib:11} p.469) but on large deviation tools.

An important idea of Ibragimov and Has'minskii has been to study the likelihood ratio $$Z_{n,\theta}(u)=\prod\limits_{i=1}^{n}\frac{f(X_i,\theta+\frac{u}{\sqrt{n}}) }{f(X_i,\theta)}~~ with~~u\in U_{n,\theta}=\{u\in\mathbb{R}~|~\theta+\frac{u}{\sqrt{n}}\in \Theta  \}.$$ Its asymptotic behavior is linked to that of the maximum likelihood estimator by the following inequality $$\mathbb{P}(\sqrt{n}(\widehat{\theta_n}-\theta)>H\mid V=\theta)\leq \mathbb{P}(\underset{|u|>H}{sup} Z_{n,\theta}(u)\geq 1\mid V=\theta) .$$
 Furthermore this quantity is connected to the Hellinger's distance:  $$\mathbb{E}[Z_{n,\theta}^{\frac{1}{2}}(u)\mid V=\theta]= 1-\frac{1}{2} r( P^{ n}_{\theta+\frac{u}{\sqrt{n}}},P^{ n}_{\theta})$$where,  for a given parametric model  $(P_{\theta})$, the  Hellinger's distance  $r$ is defined by

$$r(P_\theta,P_{\theta'})= \int (   \sqrt{f(x,\theta)} -\sqrt{f(x,\theta')}~  )^2 d\mu(x).$$  It is a measurement of the identifiability i.e the capacity  of a model to distinguish  two different values of the parameter     $\theta$.\\

The following theorem gives sufficient conditions for the consistence of the maximum likelihood estimator.\\
\begin{theorem} $\mbox{\rm{(\cite{bib:10} p.42)}}$ Let us suppose that
\begin{itemize}
\item[1)]  $\forall \theta$, $\forall n$,  the function  $u\rightarrow Z_{n,\theta}(u)$ is continuous \item[2)]  $ \forall \theta$, $\exists M>0$, $\exists m >0$ such that $\forall n$
$$ \underset{|u_1|\leq R,  |u_2|\leq R}{sup } ~|u_1-u_2|^{-{2}}~\mathbb{E}[~|Z_{n,\theta}^{\frac{1}{2}}(u_1)-Z_{n,\theta}^{\frac{1}{2}}(u_2)|^2\mid V=\theta]\leq M(1+R^m)$$
\item[3)] $\exists a>0$ such that  $\forall u\in U_{n,\theta}$,  $\forall n$ $$\mathbb{E} [Z_{n,\theta}^{\frac{1}{2}}(u)\mid V=\theta]\leq e^{-a|u|^2}.$$
Then, $\forall \theta$ , $\exists B>0$, $\exists b>0$  such that $\forall \varepsilon>0 $,  for n sufficiently large, one has  $$\mathbb{E}_{\rho}[\mathbf{1}_{\sqrt{n}|\widehat{\theta_n}(X_1,\ldots,X_n)-\theta|>\varepsilon}\mid V=\theta]\leq Be^{-b\varepsilon^2}.$$
Consequently we obtain the almost sure convergence of $\widehat{\theta_n}$ toward $\theta$.
\end{itemize}

\end{theorem} ~

\textbf{Remarks I:}  \textbf{i)} We can notice that  hypothesis  $3)$ implies the identifiability of the model. This condition is necessary because one can't find consistent estimators for a non-identifiable model.

\textbf{ii)} There exists a uniform extension  of the preceding theorem: If  $K$ is a compact set included  in $\Theta$ and if hypotheses $2)$ and $3)$ are fulfilled uniformly for $\theta\in K$ then $\exists b(K)>0$, $\exists B(K)>0$ such that $\forall \varepsilon>0 $,  for large n, $$\underset {\theta \in K}{sup}~\mathbb{E}_{\rho}[\mathbf{1}_{\sqrt{n}|\widehat{\theta_n}-\theta|>\varepsilon}\mid V=\theta]\leq Be^{-b\varepsilon^2}. \Box$$

The hypotheses of theorem $4.3.1$  may appear  restrictive, but the following result shows that they are satisfied for a large class of regular models.\\

\begin{e-proposition} $\mbox{\rm{(\cite{bib:10} p.81)}}$ If $P_\theta$ is a regular model fulfilling
\begin{itemize}
\item[1)] $0<\underset{\Theta}{inf}~J(\theta)\leq\underset{\Theta}{sup}~J(\theta)<\infty$
\item[2)] $\forall \theta$, $\forall \delta >0$ $$\underset {u\in U_{1,\theta}, |u|>\delta}{inf}{r(P_{\theta},P_{\theta+u})}>0 $$
\end{itemize}

then the hypotheses of theorem $4.3.1$ hold.

\end{e-proposition}~

From a practical point of view, the condition of local asymptotic normality introduced in the following theorem, yields a useful result for constructing confidence intervals. It possesses also a uniform version.\\
\begin{theorem} $\mbox{\rm{(\cite{bib:10} p.185)}}$ We suppose that the hypotheses of theorem  $4.3.1$ are fulfilled. Moreover  we assume that the model satisfies the local asymptotic normality condition introduced by  Le Cam: for all  $\theta$,  the sequence of stochastic processes  $(Z_{n,\theta}(u))$ converges in the sense of finite dimensional marginals toward the process $Z_{\theta}(u)=e^{u\bigtriangleup -\frac{1}{2}J(\theta)u^2}$ where $\bigtriangleup$ is a random variable distributed as $\mathcal{N}(0,J(\theta))$. Then, $\forall \theta \in \Theta$, one has

1) given $V=\theta$ $$\sqrt{n}(\widehat{\theta_n}(X_1,\ldots,X_n)-\theta)\underset{\mathcal{L}(\mathbb{P})}{\rightarrow} \mathcal{N}(0,\frac{1}{J(\theta)}),$$

2) $\forall p>0$ $$\mathbb{E}[n^{\frac{p}{2}}(\widehat{\theta_n}(X_1,\ldots,X_n)-\theta)^p\mid V=\theta]\rightarrow m_p,$$where $m_p$ is the  p-th moment of the law $\mathcal{N}(0,\frac{1}{J(\theta)})$.
 \end{theorem}~

\textbf{Remarks J:}  \textbf{i)} The hypotheses of theorem $4.3.1$ lead to the tightness of the process $(Z_{n,\theta}(u))$ in the space of  continuous functions vanishing at infinity. The pointwise convergence of this sequence becomes functional and gives $1)$.

\textbf{ii)} The maximum likelihood estimator is asymptotically unbiased and achieves asymptotically the bound of the Cramer Rao inequality.

\textbf{iii)} Since  $J$ is continuous under the hypotheses of regular models, the construction of asymptotic confidence intervals is done classically.\hfill$\Box$ \\

Now,  one of the most important property of regular models is the following:\\
 \begin{e-proposition} $\mbox{\rm{(\cite{bib:10} p.114)}}$ The condition of local asymptotic normality is fulfilled for regular models.
\end{e-proposition}~

In the following section we used those asymptotic results to give a new interpretation of $\Gamma^V_{\psi}$.

\subsubsection{\textbf{$\Gamma^{V}_{\psi}$ as an asymptotic variance}}~

 We are able to exhibit a consistent estimator in the problem of the direct estimation of $\psi(\theta)$ using the experiment generated by the observation of  $(X_1,\ldots,X_n)$, given $\psi(V)=a$. The quantity   $\Gamma_{\psi}^{V}$  will appear in the limit theorems associated to this statistical procedure.\\

\begin{e-proposition}  Under the hypotheses of proposition  $4.3.2$ one has $\forall a\in \psi(\Theta)$,

1) $\forall\varepsilon>0,$ $$\mathbb{E}[\mathbf{1}_{|\psi(\widehat{\theta_n}(X_1,\ldots,X_n))-\psi(V)|>\varepsilon}\mid \psi(V)=a] \rightarrow 0.$$

2) Given $\psi(V)=a$ $$\sqrt{n}(\psi(\widehat{\theta_n})-a)\underset{\mathcal{L}(\mathbb{P})}{\rightarrow} G_a$$
where $G_a$ is a random variable with the following density
$$g(x,a)=\mathbb{E}_{\rho}[\mathbf{1}_{\psi'\not=0}\frac{1}{\sqrt{2\pi\frac{\psi'^2}{J}}} e^{\frac{-x^2 J}{2\psi'^2}}\mid \psi=a]$$ with respect to Lebesgue measure on $\mathbb{R}$ ( $G_a$ has a variance equal to  $\Gamma^{V}_{\psi}[Id](a)$).\\
\end{e-proposition}

\textbf{Proof:} 1) We denote by  C the Lipschitz  constant of  $\psi$.

According to  Fubini theorem and by definition of the conditional expectation, $\forall (x_1,\ldots,x_n)\in (E)^n, \forall a \in \psi(\Theta)$,  $$\mathbb{E}[\mathbf{1}_{|\psi(\widehat{\theta_n}(X_1,\ldots,X_n))-\psi(V)|>\varepsilon}\mid \psi(V)=a] $$   is equal to  $$\mathbb{E}_{\rho}[\int \mathbf{1}_{|\psi(\widehat{\theta_n}(x_1,\ldots,x_n))-a|>\varepsilon} f(x_1,.)\ldots f(x_n,.)d\mu(x_1)\ldots d\mu(x_n)\mid \psi =a].$$ But we have  $$\mathbf{1}_{|\psi(\widehat{\theta_n}(X_1,\ldots,X_n))-\psi(\Theta)|>\varepsilon}\leq\mathbf{1}_{|\widehat{\theta_n}(X_1,\ldots,X_n)-\Theta|>\frac{\varepsilon}{C}} $$ and the result follows by theorem $4.3.1$ and dominated convergence theorem.\\
2)  When $\psi'(\theta_0)=0$, theorem $4.3.3$ yields   $$\mathbb{E}(\mathbf{1}_{\sqrt{n}|\psi(\widehat{\theta_n}(X_1,\ldots,X_n))-\psi(\Theta)|>\varepsilon}\mid V=\theta_0)_{\overrightarrow{n \rightarrow \infty}}~0$$ and when  $\psi'(\theta_0)\not=0$,  Slutsky's lemma (see \cite{bib:11} p.86) gives that, given  $V=\theta_0$, $$\sqrt{n}(\psi(\widehat{\theta_n}(X_1,\ldots,X_n))-\psi(\theta_0))\underset{\mathcal{L}(\mathbb{P})}{\rightarrow} \mathcal{N}\left(0,\frac{\psi'^2(\theta_0)}{J(\theta_0)}\right).$$ If  F is a bounded continuous function, using the  same argument as in  1), one has that
$$\int F(\sqrt{n}(\psi(\widehat{\theta_n})-a))\mathbb{E}_{\rho}[f(x_1,.)\ldots f(x_n,.)\mid \psi =a]d\mu(x_1)\ldots d\mu(x_n)$$ is equal to $$\mathbb{E}_{\rho}[\int F(\sqrt{n}(\psi(\widehat{\theta_n}(x_1,\ldots,x_n))-\psi))f(x_1,.)\ldots f(x_n,.)d\mu(x_1)\ldots d\mu(x_n)\mid \psi=a]$$ and the result comes by dominated convergence.$\Box$\\

\textbf{Remarks K:} \textbf{i)} When $\psi$ is injective,
$\psi(\widehat{\theta_n})$ is the maximum likelihood estimator associated to the model $(Q_a)_{a\in\psi(\Theta)}$.

\textbf{ii)}  Using the Borel-Cantelli theorem and the fact that  $\psi$ is in  $Lip^1(\Theta)$, we can extend the convergence in probability in 1) to an almost sure convergence.

\textbf{iii)}   $\Gamma^{V}_{\psi}[Id]$ is a mean of the inverse of the local Fisher information. Let us simply show this on an example: we suppose that  $\Theta=]-1,1[~\setminus~   \{0\}$, $\rho(\theta)=q(\theta)d\theta$ and  $\psi(\theta)=\theta^2$.

If  $a_0\in ]0,1[$, this point has two antecedents for $\psi$: $\theta_1=\sqrt{a_0}$  with the local Fisher information $J^{\psi(V)}_{\theta_1}(a_0)=\frac{J(\theta_1)}{\psi'^2(\theta_1)}     $ and $\theta_2=-\sqrt{a_0}$ with  $J^{\psi(V)}_{\theta_2}(a_0)=\frac{J(\theta_2)}    {\psi'^2(\theta_2)}      $. A calculus of conditional expectation gives
$$\Gamma^{V}_{\psi}[Id](a_0)=\frac{\frac{q(\theta_1))}{J^{\psi(V)}_{\theta_1}(a_0)}+\frac{q(\theta_2)}{J^{\psi(V)}_{\theta_2}(a_0)}}{q(\theta_1)+ q(\theta_2)}$$which is none other than a barycenter weighted by $\rho$.

 \textbf{iv)} When  $\rho(\theta)= q(\theta) d\theta$ with  $q$ continuous, we have similar results   if we replace the maximum likelihood estimator by the bayesian estimator associated to the quadratic loss function and the \textit{a priori} law $\rho$.

\textbf{v)} The estimation bound given in  $4.2$ becomes an asymptotic equality.\hfill$\Box$\\

In order to obtain a quadratic convergence for   $\sqrt{n} (\psi(\widehat{\theta_n})-a)$ ( allowing to  approximate in this way  $\Gamma^{V}_{\psi}[Id]$ by Monte-Carlo methods) we have to reinforce the hypotheses of proposition  $4.3.5$.\\

\begin{e-proposition} Let us suppose that the model  $(P_\theta)_{\theta\in\Theta}$ can be extended in a regular model on an open set  $\Theta'$ such that $\overline{\Theta}\subset\Theta'$. Moreover, if  \\
1) $0<\underset{\Theta'}{inf}~J(\theta)\leq\underset{\Theta'}{sup}~J(\theta)<\infty$      \\
2) $\forall \delta >0$ $$\underset{\theta\in{\Theta'}}{inf} ~\underset {u\in \tilde{U}_{1,\theta}, |u|>\delta}{inf}{r(P_{\theta},P_{\theta+u})}>0 $$where $\tilde{U}_{1,\theta}=\{u\in\mathbb{R}~|~\theta+u\in \Theta'  \}$,
then, $\forall a \in \psi(\Theta) $ $$\mathbb{E}[{n}(\psi(\widehat{\theta_n}(X_1,\ldots,X_n))-a)^2\mid \psi(V)=a]\rightarrow  \Gamma^{V}_{\psi}[Id](a).$$
\end{e-proposition}
\textbf{Proof:} Conditions 1) and 2) lead to  an uniform version of theorem  $4.3.3$: \begin{equation}\underset{\theta\in\Theta}{sup} ~(\mathbb{E}[{n}(\widehat{\theta_n}-\theta)^2\mid V=\theta]- \frac{1}{J(\theta)})\rightarrow 0.\end{equation} By Fubini theorem,  $$\mathbb{E}[{n}(\psi(\widehat{\theta_n})-a)^2\mid \psi(V)=a]$$is equal to $$\mathbb{E}_{\rho}[\int n(\psi(\widehat{\theta_n}(x_1,\ldots,x_n))-\psi)^2f(x_1,.)\ldots f(x_n,.)d\mu(x_1)\ldots d\mu(x_n)\mid \psi=a].$$ Since $\psi$ is lipschitzian, it follows from  (8) that  \begin{eqnarray*}A=\int n(\psi(\widehat{\theta_n}(x_1,\ldots,x_n))-\psi(\theta))^2f(x_1,\theta)\ldots f(x_n,\theta)d\mu(x_1)\ldots d\mu(x_n) \end{eqnarray*} fulfills $A\leq \frac{k}{J(\theta)}$ with $J^{-1}\in L^1(\rho)$ and $k\in\mathbb{R}^*_+ .$  We conclude thanks to  the dominated convergence theorem using that  $$\mathbb{E}[{n}(\widehat{\theta_n}-\theta)^2\mid V=\theta]\rightarrow \frac{1}{J(\theta)}$$ implies  
$$\hspace{5cm}\mathbb{E}[{n}(\psi(\widehat{\theta_n})-\psi(\theta))^2\mid V=\theta]\rightarrow \frac{\psi'^2(\theta)}{J(\theta)}. \hspace{5cm}\Box $$

\subsubsection{\textbf{Comments and perspectives.}}~

From the hypotheses made on the model $(P_\theta)_{\theta\in\Theta}$, we are able to give a bound concerning  the direct estimation of   $\psi(\theta)$, using the experiment generated by the observation of  $(X_1,\ldots,X_n)$ given $\psi(V)=a$. A question naturally arises: what happens when the model   $(Q_a)_{a\in\psi(\Theta)}$ is sufficiently regular to define its Fisher information matrix $J^{\psi(V)}$? One has another estimation bound that appears in some limits theorems associated to the estimation of $a=\psi(\theta)$ by means of  a n-sample of $Q_a$.

When $\psi$ is injective it is easy to show that those bounds coincide, but it is not generally the case as we can see on the following example.\\  Suppose

\begin{itemize}
\item[-] $\Theta=]-1;1[~\setminus \{0\}$$, \rho $ is distributed as the normalized uniform law on  $\Theta$
\item[-] $dP_\theta(x)=f(x,\theta)d\mu(x)= \frac{1}{\sqrt{2\pi}}e^{\frac{-(x-\theta)^2}{2}} dx$
\item[-] $\psi(\theta)=\theta^2$.\end{itemize}

The model $(P_\theta)_{\theta\in\Theta}$ is regular and fulfills the assumptions of proposition $2.2.2$ a). From the definition of the conditional expectation, we obtain for  $a \in]0,1[$
$$dQ_a(x)=\frac{f(x,\sqrt{a})+f(x,-\sqrt{a})}{2}~dx=h(x,a)dx.$$As the function $a\rightarrow h(x,a)$ is in  $C^1(]0,1[,\mathbb{R})$ and that,  according to the dominated convergence theorem, $a\rightarrow \int \frac{(h_a'(x,a))^2}{h(x,a)}dx$ is continuous, using the method of [13] p.$95$, one shows that the model $(Q_a)$ is regular.  Moreover we have $$\Gamma^{V}_{\psi}[Id](a)=4a.$$
 In order to compare  $\Gamma^{V}_{\psi}[Id]$ and $J^{\psi(V)}$ we need the following lemma.\\
\begin{lemma} Suppose  that $p(x,\theta)d\mu (x)$ and $r(x,\theta)d\mu(x)$ are two regular models on $\Theta$ such that the function $\theta\rightarrow (p(x,\theta),r(x,\theta))$ is differentiable. If we put  $s(x,\theta)= p(x,\theta)+ r(x,\theta)$ then \begin{equation}\int\frac{s'^2}{s}d\mu(x) \leq \int\frac{p'^2}{p}d\mu(x) + \int\frac{r'^2}{r}d\mu(x) .\end{equation}
\end{lemma}

\textbf{Proof:}
We set  $\tilde{s}=\sqrt{s}$, $\tilde{p}=\sqrt{p}$, $\tilde{r}=\sqrt{r}$, inequality (9) becomes \begin{equation}   \int \tilde{s}'^2 d\mu(x)\leq \int \tilde{p}'^2 d\mu(x)+\int \tilde{r}'^2 d\mu(x).\end{equation} It is easy to show that $$(10)\Leftrightarrow \tilde{p}'^2  \tilde{r}^2 + \tilde{r}'^2 \tilde{p}^2  \geq 2\tilde{p}\tilde{p}'\tilde{r}'\tilde{r} ~~\mu-a.e $$ and  $(9)$ follows with equality if and only if  $\tilde{p}\tilde{r}'=\tilde{p}'\tilde{r}$.$\Box$\\

Thus we have \begin{equation}\frac{1}{J^{\psi(V)}}>\Gamma^{V}_{\psi}[Id].\end{equation}

Hence, in this situation, we can see that error calculus gives a more precise bound. At present, we are not able to exhibit an example where $(11)$ is contradicted.
\section{\textbf{Product structures.}}

First of all, we recall the definition of the product of two error structures (see \cite{bib:2} p.200).\\

\begin{e-definition}
If  $S_i=(W_i, \mathcal{W}_i, m_i, {\mathbb{D}_i},\Gamma_i)$
(i=1,2) are two error structures, the product, denoted by $S_1 \otimes S_2$, is define as the structure
$(W_1\times W_2, \mathcal{W}_1\otimes\mathcal{W}_2,
m_1\otimes m_2, \mathbb{D}, \Gamma)$ with
\begin{center}~~~~\begin{eqnarray*}
\mathbb{D}=\{f \in L^2(m_1\otimes m_2)|~for~~m_2-almost~~every~~ y~~f(.,y)\in \mathbb{D}_1 \\
 for~~m_1-almost~~ every~~ x~~f(x,.)\in \mathbb{D}_2\\
\int \Gamma[f](x,y)dm_1(x)dm_2(y)<\infty\}\end{eqnarray*}\end{center} and
 $$\Gamma[f](x,y)=\Gamma_1[f(.,y)](x)+ \Gamma_2[f(x,.)](y).$$
\end{e-definition}~

Here we are interested in the evaluation of a  parameter   $\theta=(\theta_1,\theta_2)$ where $\theta_1$ and $\theta_2$ are supposed to be independent i.e. $V_1$ and $V_2$ are independent random variables. Let us denote by $V=(V_1,V_2):(\Omega,\mathcal{A},\mathbb{P})\rightarrow \Theta_1\times\Theta_2$ the realization of the parameter $\theta$. The law of the random variables $V$, denoted by $\rho$, fulfills:$$d\rho(\theta)=d\rho_{1}(\theta_1)d\rho_{2}(\theta_2).$$

To estimate  $\theta_1$ [resp.$\theta_2$] we choose the following regular parametric model: $$d{P}_{\theta_1}=f(x,\theta_1)~d\mu(x)~~[resp.~d{Q}_{\theta_2}=g(y,\theta_2)~ d\nu(y) ]$$ with a regular Fisher information matrix  $J_1(\theta_1)$ [resp.$J_2(\theta_2)$] such that hypothesis \textbf{E} is fulfilled.

Let us consider a random variable $X$ [resp. $Z$]  with a conditional law given $V_1=\theta_1$ [resp. $V_2=\theta_2$]
having the density: $$f(x,\theta_1)~ d\mu(x)~~~~
[resp.~g(y,\theta_2)~d\nu(y)].$$

We suppose that \textbf{$(X,V_1)$ and $(Z,V_2)$} are independent.\\

\textbf{Remark L:} We are in the situation  where the pairs (parameter, observation) are independent.   In terms of errors, this independence has to be linked with $(5)$ which is the intuitive meaning of the preceding definition of product structures.\hfill$\Box$\\

 To estimate  $\theta $, it is natural to use the conditional law of $(X,Z)$ given  $V=(\theta_1,\theta_2)$ denoted by $R_{\theta_1,\theta_2}$. From these hypotheses, it comes that
$$dR_{\theta_1,\theta_2}=f(x,\theta_1)~g(y,\theta_2)~d\mu(x)~d\nu(y).$$

Thus,  we obtain for this model the following Fisher information matrix$$
\begin{bmatrix}
J_1(\theta_1)&0\\
0&J_2(\theta_2)\\
\end{bmatrix}
$$and for $F\in Lip^1(\Theta_1\times\Theta_2),   $
$$\Gamma^{V}[F](\theta_1,\theta_2)=\frac{[F_{1}'(\theta_1,\theta_2)]^2}{J_1(\theta_1)}+
\frac{[F_{2}'(\theta_1,\theta_2)]^2}{J_2(\theta_2)}.$$

Then we have the following proposition:\\
\begin{e-proposition}

 1) $S^{V}=S^{V_1} \otimes S^{V_2}$

2) ~If $\psi_1$ and $\psi_2$ are regular changes of variables then $$(\psi_1,\psi_2)_*S^{(V_1,V_2)}=\psi_1~
_*S^{V_1}\otimes \psi_2~_*S^{V_2}.$$
\end{e-proposition}~

\noindent \textbf{Proof:} 1) Let us notice that  $Lip^1(\Theta_1\times\Theta_2)$ is included in the domain of the product structure $S^{V_1} \otimes S^{V_2}$.

Moreover, from the expression of the information matrix, for $F\in Lip^1(\Theta_1\times\Theta_2)  $ it follows that $$  {\Large \mathcal{E}}^{V}   [F]=
\int {\Large
\mathcal{E}}^{V_1}[F(.,y)]d\rho_2(y)+\int {\Large
\mathcal{E}}^{V_2}[F(x,.)]d\rho_1(x).$$ Thus  $\parallel . \parallel_{\mathcal{E}^{V}}$ coincides on $Lip^1(\Theta_1\times\Theta_2)$ with the norm associated to the product structure. Hence, we can deduce that the hypothesis \textbf{E} is fulfilled for the model $(R_{\theta_1,\theta_2})_{\Theta_1\times\Theta_2}$: $S^V$ is well-defined.

Furthermore, since $Lip^1(\Theta_1\times\Theta_2)$ is dense in $(\mathbb{D}^{V}, \parallel . \parallel_{\mathcal{E}^{V}} )$,    $\mathbb{D}^{V}$ is included in the domain of the product structure and the two squared field operators coincide on $\mathbb{D}^{V}$.

For the other inclusion, we use the fact that the functions of the form  $F=\sum\limits_{i=1}^p f_i g_i$ with $f_i\in\mathbb{D}^{V_1}$ and $g_i \in\mathbb{D}^{V_2}$ are dense in the domain of the product structure for the associated norm (see \cite{bib:2} p.201) and belong to $\mathbb{D}^V$ as easily seen using the closedness of the forms $\mathcal{E}^{V_1}$ and $\mathcal{E}^{V_2}$.

2) The equality comes from 1) and section 3.\hfill$\Box$\\

\textbf{Remarks M:} \textbf{i)} The preceding results  extended obviously to n-tuple.

\textbf{ii)} We can notice that this property expresses the additive property of the Fisher information for independent experiments.\hfill$\Box$ \\

Since it is easy to build infinite products of error structures  (see \cite{bib:2},\cite{bib:5},\cite{bib:**}), we are able to obtain an  empirical error calculus associated to the estimation of the parameters of the type    $\theta=(\theta_i)_{i\in\mathbb{N}}$ working component per component.

\section{The choice of an \textit{a priori} law $\rho$.}

In the preceding sections, the choice of an \textit{a priori} law on the  space of parameters $\Theta$ is left to the  practitioner as in the bayesian analysis. The determination of our error structure $S^{V}$ can appear, to some degree, incomplete. We are going to show that, once a regular parametric model is chosen, a natural probability measure becomes apparent: the Jeffreys prior (see \cite{bib:11} p.490).  This probability is well known in bayesian analysis. Moreover, it possesses a remarkable stability concerning error calculus: it is  invariant under reparameterization  and compatible with the notion of product.

Let  $(P_\theta)_{\theta\in\Theta}$ be a regular model such that       $$ K=\int  _{\Theta} \sqrt{det(J(\theta))} d\theta < \infty .$$  We can define on $\Theta$ the following probability measure $$\rho^{V}(d_\theta)=\frac{\mathbb{I}_\Theta \sqrt{det(J(\theta)) }d\theta}{K} $$called the Jeffreys prior induced by   the model $(P_\theta)_{\theta\in\Theta}$. It is often used in bayesian analysis for its invariance under reparameterisation. Moreover it is the prior measure which has the smallest influence on the posterior measure in the sense of the asymptotic Shannon information (see \cite{bib:12}). In term of error calculus its properties are summarized in the following proposition:\\
\begin{propo}a) If $\psi:\Theta\rightarrow \mathbb{R}^d$ is a regular change of variables $$\psi_*\rho^V=\rho^{\psi(V)}.$$b) In the framework of section 5 $$\rho^{(V_1,V_2)}=\rho^{V_1}\otimes \rho^{V_2}.$$
\end{propo}
\noindent \textbf{Proof:} Obvious using the classical properties of the Fisher information.\hfill$\Box$\\

 %\textbf{Remark N:} When d=1 and $\rho=\rho^V$,  Hamza theorem ensure that $\textbf{E}$ is always fulfilled. In general, if  %$\theta\rightarrow  J(\theta)$ is differentiable, we can show, using an integration by parts, the existence of an operator  %$S$  from $Lip^1(\Theta)$ into  $L^2(\mu^V)$ such that $\mathcal{E}^V(F,G)=<S[F],G>_{L^2(\mu^V)}$. Thus $\textbf{E}$ holds %by a classical argument (see).$\bullet$

Finally, with suitable hypotheses, the Jeffreys prior may be seen as the invariant measure of the generator associated to the induced infinitesimal perturbation in the convergence of the maximum likelihood estimator.

\section{\textbf{Conclusion.}}

Through statistical experiments, we have seen that the fundamental identification gave an error structure intrinsically linked to the observed physical phenomenon. The remarkable robustness of this identification, regarding injective changes of variables and  products, yields a particularly efficient tool for finite dimensional estimation.

The existence of such an error structure built from the parametric model allows to propagate the accuracy through calculations performed with the parameter thanks to a coherent specific differential calculus (property $1$ of $\Gamma$). Moreover error calculus provides a natural framework  concerning the study of non-injective mapping.  A possible extension   will be to generalize such an experimental protocol when $J$ is singular and also to explore more precisely the connections between Dirichlet forms and  asymptotic statistics. Finally, we wonder whether   the semi-parametric and non-parametric estimation theories (see \cite{bib:13}) could lay the foundation of an infinite dimensional identification in order to get $\Gamma$ on the Wiener space, using a direct functional reasoning instead of  a component per component argument as above.

\end{document}